\documentclass[12pt]{amsart} 
\usepackage{lmodern}
\usepackage{amsmath}
\usepackage{amssymb}
\usepackage{latexsym}

\def\Nr{{\mathbb N}}

\def\Rr{{\mathbb R}}

\def\Pc{{\mathcal{P}}}

\def\Sc{{\mathcal{S}}}

\def\one{{\rm \bf 1}}

\def\dist{\operatorname{dist}}

\def\conv{\operatorname{conv}}
\def\({\left(}     
\def\){\right)}    
\def\[{\left[}     
\def\]{\right]}

\newtheorem{definition}{Definition}
\newtheorem{theorem}{Theorem}
\newtheorem{lemma}{Lemma}
\newtheorem{proposition}{Proposition}
\newtheorem{remark}{Remark}

\begin{document}
\author{Freddy Delbaen}
\address{Departement f\"ur Mathematik, ETH Z\"urich, R\"{a}mistrasse
101, 8092 Z\"{u}rich, Switzerland, also Institut f\"ur Mathematik,
Universit\"at Z\"urich, Winterthurerstrasse 190,
8057 Z\"urich, Switzerland}
\email{delbaen@math.ethz.ch}
\title{Convex Increasing Functionals on $C_b(X)$ Spaces}

\maketitle

\abstract  We prove that convex functions on a $C_b(X)$ space satisfying a mild continuity condition  can be represented using sigma additive measures. This generalises a result of Cheridito, Kupper and Tangpi, \cite{CKT}.\endabstract

\section{Notation and Preliminaries}

We first recall some definitions from general topology.  A Polish space $X$, is a topological space that is separable and where the topology is metrisable with a complete metric. In this paper we will always use $X$ for a fixed Polish space. No extra assumptions will be put on $X$.  There are no local compactness assumptions as we believe that the most interesting applications are when $X$ is for instance an infinite dimensional Banach space.  The space of real valued bounded continuous functions defined on $X$ is denoted by $C_b(X)$, the space of all continuous functions is denoted $C(X)$. The aim of the paper is to find a representation theorem for convex, monotonic functions defined on $C_b(X)$. We say that a function $F\colon C_b(X)\rightarrow \Rr$ is nondecreasing if for elements in $C_b(X)$, $f\le g$ we always have $F(f)\le F(g)$. Such an assumption immediately shows that $F$ is continuous for the sup norm on $C_b(X)$.  Indeed if $\Vert f-g\Vert\le\epsilon$ (where $\epsilon>0$ and $\Vert .\Vert$ denotes the sup norm), then monotonicity implies 
\begin{eqnarray*}
F(f-\epsilon)&\le & F(g)\le F(f+\epsilon)\\
F(f-\epsilon)&\le & F(f)\le F(f+\epsilon)\\
\text{ hence }& & |F(f)-F(g)|\le F(f+\epsilon)-F(f-\epsilon)
\end{eqnarray*}
Since the function $\Rr\rightarrow \Rr; a\rightarrow F(f+a)$ is convex, it is continuous, so we find $\lim_{\epsilon\rightarrow 0}F(f+\epsilon)-F(f-\epsilon)=0$.
\begin{definition} We will say that a nondecreasing function $F\colon C_b(X)\rightarrow \Rr$ is upward continuous if $f_n\uparrow f$ in $C_b(X)$, pointwise,  implies $F(f_n)\uparrow F(f)$.
\end{definition}
Cheridito, Kupper and Tangpi, \cite{CKT} have studied  convex nondecreasing functions that were downward continuous. They showed that such  functions can be represented by sigma--additive measures on $X$. Since downward continuity implies upward continuity, our setup is more general. To prove this let $h_n\downarrow 0$ (pointwise) be a sequence in $C_b(X)$ and take $f\in C_b(X)$.  Convexity and monotonicity imply that
$$
0\le F(f)-F(f-h_n)\le F(f+h_n)-F(f),
$$
from which it follows that downward continuity implies upward continuity.  In this paper $F$ will always be a real valued nondecreasing upward continuous convex function defined on the whole space $C_b(X)$. To simplify the statements, we will also assume that $F(0)=0$.  This is no restriction as we can always replace $F$ by a the new function $f\rightarrow F(f)-F(0)$. It avoids that the value taken at $0$ always shows up and complicates the writing.  The present author is known to be lazy.

The space $M(X)$ is not the dual space of $C_b(X)$, the weak topology $\sigma(C_b(X),M(X))$ and the Mackey topology  $\tau(C_b(X),M(X))$ are not metrisable. Since the requirement for the convex functions will be stated for sequences, the non-metrisability poses a first problem. The duality theory of convex functions, see \cite{phelps}, cannot be used in our case.  On the space $M(X)$ we will use the weak$^*$ topology $\sigma(M(X),C_b(X))$, see \cite{Par} for details where it is called weak topology as is common in the probabilistic literature\footnote{We prefer to call it weak$^*$ since it is the restriction of $\sigma(C_b(X)^*,C_b(X))$ to $M(X)$.  In functional analysis the weak topology refers to the finer topology $\sigma(M(X),M(X)^*)$.}.

The present paper is a generalisation of the author's paper on monetary utility functions, \cite{C(X)}, satisfying (because of concavity) a downward continuity assumption. The author is grateful to Patrick Cheridito and Michael Kupper for insisting on looking at the more general context of monotonic convex functions. We are also grateful to Matteo Burzoni for mentioning related papers.
\section{Some topological Results}

We will not use the dual space of $C_b(X)$ which involves measures on the Stone-\v Cech compactification of $X$. Instead we will use another more constructive approach. It is known that a Polish space can be embedded as a $G_\delta$ in a compact metric space $K$\footnote{We could have used the Stone-\v Cech  compactification, $\beta X$, of $X$. But we used a more constructive approach to please the fundamentalists.}. These results were proved by Alexandroff, Hausdorff and Sierpinsky, see the book by Kelley, \cite{Kelley}, p. 208, Problem K. The set $K\setminus X$ is the union of a sequence of compact sets $(L_n)_n$. 
\begin{lemma}\label{phi-n} For each set $L_n$ there is an element $\phi_n\in C(K)$ such that $0\le \phi_n\le 1$  and $L_n=\phi_n^{-1}(1)$.  In particular $\phi_n(x)<1$ for each $x\in X$. 
\end{lemma}
{\bf Proof } We only have to give a construction of the functions $\phi_n$. On $K$ we use a compatible metric, $d$, and the distance function to $L_n$, $\dist(k,L_n)=\min \{d(k,l)\mid l\in L_n\}$.  Then we define the function $\phi_n(k) = \max(1- \dist(k,L_n),0)$.\qed

Not every function in $C_b(X)$ can be extended to a continuous function on $K$. If $d$ is a metric defining the topology on $K$, a function in $C_b(X)$ can be extended to a continuous function on $K$ if and only if it is uniformly continuous (for the metric $d$) on $X$. Clearly we have $C(K)$ is an isometric subspace of $C_b(X)$.  However we also have.
\begin{lemma} \label{approx} For each element $f\in C_b(X)$ there is a uniformly bounded sequence of functions $g_n\in C(K)$ such that $g_n\uparrow f$ on $X$.
\end{lemma}
{\bf Proof } The proof  uses a convolution technique.  For each $n$ and each $x\in K$, we define
$$
g_n(x)=\inf \{f(y)+n\,d(x,y)\mid y\in X\}.
$$
Clearly if $x\in X$, $-\Vert f\Vert \le g_n(x)\le f(x)\le \Vert f\Vert$.  For two points $x_1,x_2$ in $K$ we have
\begin{eqnarray*}
g_n(x_1)&=&\inf\{f(y)+n\,d(x_1,y)\mid y\in X\}\\
&\le& \inf \{f(y)+n\,d(x_2,y) + n\,d(x_1,x_2)\mid y\in X\}\\
&=& g_n(x_2)+n\,d(x_1,x_2).
\end{eqnarray*}
By symmetry we then get $|g_n(x_1)-g_n(x_2)|\le n\, d(x_1,x_2)$, proving that $g_n$ is an element of $C(K)$. It is also obvious that $g_n\le g_{n+1}$.  We now show that for each $x\in X$, $g_n(x)\uparrow f(x)$.  To estimate $g_n(x)$ we distinguish between $d(x,y)\le \frac{2\Vert f \Vert}{n}$ and $d(x,y)\ge \frac{2\Vert f \Vert}{n}$.  In the latter case we have $f(y)+n\,d(x,y)\ge f(y) + 2\Vert f\Vert\ge \Vert f\Vert\ge f(x)$.  In the former case we have $f(y)+n\,d(x,y)\ge \inf_{d(x,z)\le \frac{2\Vert f \Vert}{n}}f(z)$. We get $g_n(x)\ge \inf_{d(x,y)\le \frac{2\Vert f \Vert}{n}}f(y)$. As $n\rightarrow \infty$ the continuity of the function $f$ implies that $g_n(x)\rightarrow f(x)$.\qed

The approximation theorem actually shows that to calculate $F^*(\mu)$ for $\mu\in M(X)$ we can restrict $F$ to the space $C(K)$. The restriction of $F$ to $C(K)$ is norm continuous and hence can be represented by its conjugate function defined on the dual space of $C(K)$, known as the space, $M(K)$, of sigma additive measures on the Borel sigma algebra of $K$. We define this conjugate function as
\begin{eqnarray*}
F^*(\mu)&=&\sup\{\mu(f)-F(f)\mid f\in C(K)\}\text{ and for $\mu\in M(X)$:}\\
&=&\sup\{\mu(f)-F(f)\mid f\in C_b(X)\}\text{ as we observed above}\\
F^*&\colon& M(K)\rightarrow \Rr\cup\{+\infty\}\\
\mu(g)&\le& F(g)+F^*(\mu)\text{ for }\mu\in M(X), g\in C_b(X).
\end{eqnarray*}
The classical duality, \cite{phelps}, then says for $f\in C(K)$:
$$
F(f)=\sup\{\mu(f)-F^*(\mu)\mid \mu \in M(K)\}.
$$
Of course we only need the measures $\mu$ for which $F^*(\mu)<+\infty$.  The following proposition is important in our setting
\begin{proposition} The conjugate function $F^*$ satisfies
\begin{enumerate}
\item For each $\mu\in M(K)$ we have $0\le F^*(\mu)\le +\infty$.
\item $F^*(\mu)<\infty$ implies that $\mu\ge 0$.
\item $F^*$ is convex and is lower semi continuous for the weak$^*$ topology  $\sigma(M(K),C(K))$.
\item The epigraph of $F^*$\!, $Ep=\{(\mu,t)\mid F^*(\mu)\le t <+\infty\}$ is a convex closed set in $M_+(K)\times \Rr_+$, where $M(K)$ has the weak$^*$ topology $\sigma(M(K),C(K))$. For this topology the epigraph is a Polish space.
\end{enumerate}
\end{proposition}
{\bf Proof } Except for the last sentence, this is standard but for completeness we give a sketch of the proof. Because $F(0)=0$ we immediately get $F^*(\mu)\ge 0$.  If the measure $\mu$ is not a positive measure, there is a nonpositive function $f$ such that $\mu(f)>0$.  Since $F$ is monotonic we  have $F(f)\le 0$ and hence we have $F^*(\mu)\ge \mu(f)$.  Since this also holds for $\lambda f$ for each $\lambda >0$ we get $F^*(\mu)=+\infty$. For fixed $\alpha\in\Rr$, the sublevel set $\{\mu\mid F^*(\mu)\le\alpha\}$  is clearly closed since it equals $\cap_{f\in C_b(X)}\{\mu\mid \mu(f)-F(f)\le\alpha\}$.  That implies that the epigraph is closed.  That it is a Polish space follows from the fact that the set $M_+(K)$ endowed with the weak$^*$ topology is a Polish space. 
\qed
\begin{remark}{\rm The set $M(K)$ endowed with the weak$^*$ topology is not a Polish space, it is not even metrisable,  but the cone of nonnegative elements, $M_+(K)$,  with the weak$^*$ topology forms a Polish space.
}
\end{remark}
\section{The Representation Theorem}

The formula proved using elements $\mu\in M_+(K)$, cannot be used to represent all functions of $C_b(X)$.  There are two shortcomings.  The first is that it is only valid for elements of the strictly smaller space $C(K)$, second it uses nonnegative measures that might charge the set $K\setminus X$. The proof of the main result addresses these two shortcomings.
\begin{theorem} For a nondecreasing upward continuous convex function $$F\colon C_b(X)\rightarrow \Rr,\quad F(0)=0,$$ the conjugate function, defined on $M_+(X)$, satisfies
\begin{eqnarray*}
F^*&\colon& M_+(X)\rightarrow \Rr_+\cup\{+\infty\}\\
F^*(\mu)&=&\sup\{\mu(f)-F(f)\mid f\in C(K)\}\\
&=&\sup\{\mu(f)-F(f)\mid f\in C_b(X)\},
\end{eqnarray*}
and we have for all $f\in C_b(X)$:
$$
F(f)=\sup\{\mu(f)-F^*(\mu)\mid \mu \in M_+(X)\}.
$$
\end{theorem}
{\bf Proof } The approximation Lemma, \ref{approx} shows that the calculation of $F^*$ can be done using all elements of $C_b(X)$ and not just the elements of $C(K)$. In what follows we first work with the conjugate function defined on $M_+(K)$. We now show that the convex set
$$
E=\{(\mu,t)\mid \mu\in M_+(X),+\infty> t\ge F^*(\mu)\}
$$
is weak$^*$ dense in the epigraph  $$Ep=\{(\mu,t)\mid \mu\in M_+(K), F^*(\mu)\le t <+\infty\}.$$ For each $\epsilon>0$ and each $n$ we use the convex sets
$$
W_{n,\epsilon}=\{(\mu,t)\mid (\mu,t)\in Ep, \mu(L_n)<\epsilon\},
$$
where $L_n, \phi_n$  refer to Lemma \ref{phi-n}.  The definition of the weak$^*$ topology for measures shows that $W_{n,\epsilon}$ is an open set in $Ep$.  We will show that it is weak$^*$ dense in $Ep$.  Suppose on the contrary that there is an element $(\mu_0,t_0)\in Ep\setminus \overline{W_{n,\epsilon}}$. According to the Hahn Banach theorem we can find an element $f\in C(K), \alpha\in\Rr$ such that
$$
\mu_0(f)+\alpha t_0 <\inf\{ \mu(f)+\alpha t\mid (\mu,t)\in W_{n,\epsilon}\}.
$$
In case $W_{n,\epsilon}$ is empty the right side is $+\infty$ and we can take $(f,\alpha)=(1,1)$. In case the set is nonempty, the number $\alpha\ge 0$ since otherwise the right side would $-\infty$. The inequality then changes into
$$
\mu_0(f)+\alpha F^*(\mu_0) <\inf\{ \mu(f)+\alpha F^*(\mu)\mid F^*(\mu)<+\infty, \mu(L_n)<\epsilon\}.
$$
Since $0\le F^*(\mu_0)<\infty$ we can increase $\alpha$ a little bit and still have the same inequality.  In other words we may suppose that $\alpha>0$.  Then we can divide by $\alpha$ to get a new function $g\in C(K)$ such that
$$
\mu_0(g)+F^*(\mu_0)<\inf\{ \mu(g)+ F^*(\mu) \mid F^*(\mu)<+\infty, \mu(L_n)<\epsilon\}.
$$
Changing signs we get numbers $\gamma\in\Rr, \delta>0$ such that
\begin{eqnarray*}
& &\mu_0(-g) - F^*(\mu_0)>\gamma+\delta\\
 &\ge& \gamma=\sup\{ \mu(-g)- F^*(\mu) \mid \mu\in M_+(K), \mu(L_n)<\epsilon\}.
\end{eqnarray*}
Let us now fix a number  $k$ so that $k\epsilon > \delta$. Looking at functions 
$g_m=-g-k\phi_n^m$ we get for $\mu(L_n)\ge \epsilon$ that $\mu(-g-k\phi_n^m)-F^*(\mu)=\mu(-g)-F^*(\mu)-k\mu(\phi_n^m)\le F(-g)-\delta$ whereas for $\mu(L_n)<\epsilon$ we have 
\begin{eqnarray*}
\mu(-g-k\phi_n^m)-F^*(\mu)&\le&\mu(-g)-F^*(\mu)\\
&\le&\gamma\le \mu_0(-g)-F^*(\mu_0)-\delta\\
&\le& F(-g)-\delta.
\end{eqnarray*}
Summarizing we get for all $\mu$ that $\mu(-g-k\phi_n^m)-F^*(\mu)\le F(-g) -\delta$.  Since $-g-k\phi_n^m\in C(K)$, we get that $F(-g-k\phi_n^m)\le F(-g)-\delta$. By passing to $m\rightarrow +\infty$ and using upward continuity  $g_m\uparrow -g$ for sequences converging only on $X$, we find that $F(-g)\le F(-g)-\delta$, a contradiction.  This shows the density of $W_{n,\epsilon}$. Since $Ep$ is a Polish space the countable intersection of dense $G_\delta$ sets is still dense and so we find that
$E$ is dense in $Ep$.  This immediately gives that for $f\in C(K)$ we have
$$
F(f)=\sup\{\mu(f)-F^*(\mu)\mid \mu \in M_+(X)\}.
$$
The approximation Lemma \ref{approx} and the upward continuity then also show that for all $f\in C_b(X)$:
$$
F(f)=\sup\{\mu(f)-F^*(\mu)\mid \mu \in M_+(X)\}.
$$
\qed
\begin{proposition} The function $F^*\colon M(X)\rightarrow \overline{\Rr_+}$  is lower semi continuous.  For each $m<\infty$ the set $\Sc_m=\{\mu\mid F^*(\mu)\le m\}$ is convex, weak$^*$ closed and bounded.
\end{proposition}
{\bf Proof } The first part is immediate since $F^*(\mu)=\sup\{\mu(f)-F(f)\mid f\in C_b(X)\}$. The boundedness follows from the inequality $\Vert \mu\Vert=\mu(1)\le F^*(\mu)+F(1)$. \qed

For convex nondecreasing functions on $C_b(X)$ we can extend the upward continuity in the following way.
\begin{proposition} If $(f_n)_n$ is a uniformly bounded sequence in $C_b(X)$ converging pointwise to $f\in C_b(X)$ then $F(f)\le \liminf F(f_n)$.
\end{proposition}
{\bf Proof } We use the Fatou lemma from measure theory. First we observe that for every $\epsilon>0$ we can find a nonnegative measure $\mu$ such that
$F(f)\le \mu(f)-F^*(\mu)+\epsilon$. This implies that  $F(f)\le \lim \mu(f_n)-F^*(\mu)+\epsilon$ which then implies that $F(f)\le \liminf F(f_n)+\epsilon$. Since $\epsilon$ was arbitrary this ends the proof.\qed

\section{The relation with weak$^*$ compactness}

We already mentioned that downward continuity implies upward continuity.  In this section we will show that downward continuity is related to weak$^*$ compactness properties.  We start with some general remarks on the growth of the function $F^*$. Let $F$ be a nondecreasing upward continuous convex function. Then the epigraph $E=\{(f,t)\mid F(f)\le t <+\infty\}$ is a convex set in $C_b(X)\times \Rr$.  It is closed for the product topology of the norm topology on $C_b(X)$ and the usual topology on $\Rr$. It has a nonempty interior but this is not very useful since $M(X)$ is not the dual of $C_b(X)$ for the norm topology. Nevertheless we will show that there are supporting functionals.
\begin{proposition} \label{bound}For every $1\ge\epsilon >0$ and $f\in C_b(X)$ there is an $\epsilon-$supporting functional, $\mu\in M_+(X)$, for the epigraph. That means that for every $g$:
$$
F(f+g)\ge F(f)-\epsilon +\mu(g).
$$
Furthermore we have that
$$
\mu \text{ is $\epsilon-$supporting if and only if } F(f)\le \mu(f)-F^*(\mu)+\epsilon.
$$
There also is a function $\psi\colon \Rr_+\rightarrow \Rr_+$, independent of $0<\epsilon\le 1$ such that for $\mu$, $\epsilon-$supporting, $\Vert \mu\Vert, F^*(\mu) \le \psi(\Vert f\Vert)$.
\end{proposition}
{\bf Proof } To get uniform bounds we only consider $1\ge \epsilon>0$. The relation $F(f)=\sup\{\mu(f)-F^*(\mu)\mid \mu\in M(X)\}$ shows that the epigraph is closed for the product topology of $\sigma(C_b(X),M(X))$ and the topology on $\Rr$. For given $\epsilon >0$ and $f\in C_b(X)$ the point $(f,F(f)-\epsilon)$ is not in the epigraph. The Hahn-Banach theorem then provides a separating hyperplane or in other words, there is $\mu\in M(X)$ satisfying for all $g\in C_b(X)$:
$$
F(f+g)\ge F(f)-\epsilon +\mu(g).
$$
Such a measure $\mu$ is nonnegative since otherwise we could find $g_n\le 0$ such that $\mu(g_n)\rightarrow +\infty$ whereas $F(f+g_n)$ would remain smaller than $F(f)$. We can also argue as follows:
\begin{eqnarray*}
& &\forall g: F(f+g)\ge F(f)-\epsilon +\mu(g)\\
&\Longleftrightarrow& \forall g: F(f)\le \mu(f)-(\mu(f+g)-F(f+g))+\epsilon\\
&\Longleftrightarrow& F(f)\le \mu(f)-\sup\{(\mu(f+g)-F(f+g))\mid g\in C_b(X)\}+\epsilon\\
&\Longleftrightarrow& F(f)\le \mu(f)-F^*(\mu) +\epsilon.
\end{eqnarray*}
If $\mu$ is $\epsilon-$supporting at $f$ then $\Vert\mu\Vert=\mu(1)\le F(f+1)-F(f)+\epsilon$ which by monotonicity is smaller than $F(\Vert f\Vert +1)-F(-\Vert f\Vert)+1$. Also $F^*(\mu)\le \mu(f)-F(f)+\epsilon$ yields another bound only depending on $\Vert f\Vert$. Summarising there is a function $\psi$ as required.\qed

From \cite{CKT} we recall
\begin{proposition} A  nondecreasing convex function $F\colon C_b(X)\rightarrow \Rr$, $F(0)=0$, that is downward continuous at $0$ is downward continuous at every $f$.  Consequently, a nondecreasing convex function that is downward continuous at one point is downward continuous at all points.
\end{proposition}
{\bf Proof }  Before proceeding with the proof let us first observe that the function $\Rr\rightarrow\Rr; a\rightarrow F(af)$ is convex and hence continuous. Take now a sequence, $(h_n)_n$, in $C_b(X)$ such that $h_n\downarrow 0$ pointwise and suppose that for such sequences $F(h_n)\downarrow 0$. Then take $1>\epsilon>0$ and write $f+h_n=(1-\epsilon) \frac{1}{1-\epsilon}f + \epsilon \frac{1}{\epsilon}h_n$.  Convexity 
implies
$$
F(f+h_n)\le (1-\epsilon)F\( \frac{1}{1-\epsilon}f\) + \epsilon F\(\frac{1}{\epsilon}h_n\).
$$
For $n\rightarrow 0$ downward continuity at 0 gives
$$
\lim_n F(f+h_n)\le (1-\epsilon)F\( \frac{1}{1-\epsilon}f\).
$$
This holds for every $\epsilon>0$.  Continuity with respect to $\epsilon$ then gives 
$\lim_n F(f+h_n)\le F(f)$. Conversely if $F$ is downward continuous at $f$,  we look at the function $G(g)=F(f+g)-F(f)$ which is monotonic, convex and downward continuous at $0$, hence downward continuous everywhere. For $g=-f$ this just says that $F$ is downward continuous at $0$.\qed
\begin{remark}{\rm  A similar argument can be used to prove the same statement for upward continuous functions.
}\end{remark}
We can now prove the main result of this section.
\begin{theorem} \label{compact}A nondecreasing convex function $F:C_b(X)\rightarrow \Rr$ with  $F(0)=0$ is downward continuous if and only if for each $m$ the set $\Sc_m=\{\mu\mid F^*(\mu)\le m\}$ is weak$^*$  compact.
\end{theorem}
\begin{remark}{\rm We already showed that $\Sc_m$ is convex, bounded and  weak$^*$ closed, even when $F$ is only upward continuous. The compactness property is by Prohorov's theorem, \cite{Par} theorem 6.7 page 47, equivalent to uniform tightness. In \cite{Par} only probability measures are treated but the proof for bounded sets of nonnegative measures is almost the same.
}\end{remark}
\begin{lemma} A bounded set $\Gamma\subset M_+(X)$ is uniformly tight if and only if for each decreasing sequence $h_n\downarrow 0$ in $C_b(X)$ we have
$$
\sup \{\mu(h_n)\mid \mu\in \Gamma\}\rightarrow 0.
$$
\end{lemma}
{\bf Proof of the lemma } If $\Gamma$ is uniformly tight then for all $\epsilon>0$ there is a compact set $L\subset X$ such that $\sup_{\mu\in\Gamma}\mu(L^c)\le\epsilon$.  Dini's theorem shows that when $h_n\downarrow 0$ pointwise on X,  the convergence is uniform on $L$.  Hence there is $n_0$ so that $h_n(x)\le \epsilon$ for all $x\in L$ and all $n\ge n_0$.  For such $n$ we  have for all $\mu\in \Gamma$: $\mu(h_n)\le \epsilon \mu(L)+\Vert h_1\Vert \epsilon$.  This shows that $\sup \{\mu(h_n)\mid \mu\in \Gamma\}\rightarrow 0$.  Conversely we will show that for each $\delta>0$ and each $\epsilon>0$ there is a finite number of balls of radius $\delta$ so that the complement of their union has measure smaller than $\delta$ for each measure in $\Gamma$.  According to the proofs in \cite{Par} especially the remark on page 49, this shows uniform tightness.  We now proceed with the details.  We take a metric on $X$ for which $X$ is complete.  We also take a countable dense set $\{x_k\mid k\ge 1\}$ in $X$.  For given $\delta$ the union $\cup_k B(x_k,\delta)$ of the open balls then covers $X$. For each $k$ we take a function $\phi_k\in C_b(X)$ satisfying $0\le\phi_k(x)<1$ for $x\in B(x_k,\delta)$ and $\phi_k(x)=1$ outside $B(x_k,\delta)$. For each $n$ we now define $h_n=\Pi_{k=1}^n\phi_k^n$. Clearly $h_n\downarrow 0$ on $X$. The hypothesis then implies that for $\epsilon>0$, there is $n$ such that $\sup \{\mu(h_n)\mid \mu\in \Gamma\}\le\epsilon$, Hence also $\sup \{\mu\( \(\cup_{k=1}^nB(x_x,\delta) \)^c\)\mid \mu\in \Gamma\}\le \epsilon$.\qed

{\bf Proof } We first suppose that $F$ is downward continuous and we will show that $\Sc_m$ is uniformly tight. Suppose that uniform tightness does not hold.  According to the lemma there is  a $\delta>0$, a sequence $\mu_n\in \Sc_m$, a sequence $(h_n)_n \in C_b(X)$ with $h_n\downarrow 0$ such that $\mu_n(h_n)\ge \delta$.  Take now $k$ such that $k\delta-m\ge 1$.  Then $k\,h_n\downarrow 0$ but
$$
F(k\,h_n)\ge \mu_n(k\, h_n)-F^*(\mu_n)\ge k\delta-m\ge 1,
$$
A  contradiction to the downward continuity.  Suppose now conversely that $\Sc_m$ is weak$^*$ compact. Take a sequence $h_n\downarrow 0$ in $C_b(X)$.  Take a sequence $1\ge \epsilon_n>0$ such that $\epsilon_n\rightarrow 0$. For each $n$ take $\mu_n$ such that $F(h_n)\le \mu_n(h_n)-F^*(\mu_n)+\epsilon_n$. We already know by Proposition \ref{bound} that for all $n$, $F^*(\mu_n)\le \psi(\Vert h_1\Vert)$ and hence the sequence $\mu_n$ is uniformly tight.  We then have $\mu_n(h_n)\rightarrow 0$.  This shows that $F(h_n)\rightarrow 0$ (and at the same time $F^*(\mu_n)\rightarrow 0$).\qed
\begin{theorem}  If $F$ is a nondecreasing convex function $F\colon C_b(X)\rightarrow \Rr$ with $F(0)=0$, if moreover $F$ is downward continuous then for each $f$ there is a $\mu\in M_+(X)$ such that $\mu(f)=F(f)+F^*(\mu)$.  Consequently the epigraph of $F$ has a supporting functional at each point $(f,F(f))$.
\end{theorem}
{\bf Proof }  Take $f\in C_b(X)$. For each $\epsilon>0$, there is an $\epsilon-$supporting functional, given by a measure $\mu_\epsilon$.  But from the previous theorem we know that this family of measures is uniformly tight since we have the bound $\sup_{1\ge\epsilon>0}F^*(\mu_\epsilon)\le \psi(\Vert f\Vert)<\infty$.  Hence there is a sequence $\epsilon_n\rightarrow 0$ such that the sequence $\mu_{\epsilon_n}$ converges weak$^*$ to, say, a measure $\mu$.  Clearly the inequalities
$$
F(f+g)\ge \mu_{\epsilon_n}(g) +F(f)-\epsilon_n,
$$
imply that $F(f+g)\ge \mu(g) +F(f)$, valid for all $g\in C_b(X)$. This is equivalent to $\mu(f)=F(f)+F^*(\mu)$.\qed

We will now prove the converse of the previous theorem. The proof uses a result which is called the sup-limsup theorem.  The version we need was proved by Gal\'an and Simons, \cite{gal-sim}.  Such results are used to prove the James's weak compactness theorem and variants of it. Because we are dealing with the special case of measures and continuous functions, we do not need to adapt the technique of the ``undetermined'' functions, see \cite{gal-sim} for more information and other examples of the use of this technique. Let us first recall Theorem 7 of Gal\'an and Simons.\footnote{The author wants to thank Professor Orihuela for introducing him to the work of Gal\'an and Simons.}

\begin{theorem}  Let $E$ be an arbitrary set and let $G\subset E$ be a nonempty subset. Let $(f_k)_k$ be a uniformly bounded sequence of functions $f_k\colon E\rightarrow \Rr$. By $\conv_\sigma(f_k;k\ge 1)$ we mean the set:
$$
\conv_\sigma(f_k;k\ge 1)=\left\{\sum_{k\ge 1}\lambda_k f_k\mid \lambda_k\ge 0,\sum_k\lambda_k=1\right\}.
$$
Let $G$ satisfy the following peak condition:
$$
\forall f\in \conv_\sigma(f_k;k\ge 1)\,\, \exists x\in G\text{ such that }f(x)=\sup_{y\in E} f(y).
$$
Then  the function $\limsup_k f_k$ satisfies:
$$
\sup \{\limsup_k f_k(x);x\in G\}=\sup \{\limsup_k f_k(x);x\in E\}.
$$
\end{theorem}

\begin{theorem} Let $F\colon C_b(X)\rightarrow \Rr$ be a nondecreasing convex function, $F(0)=0$, that satisfies the upward continuity condition. Suppose that for each $f\in C_b(X)$ there is $\mu\in M_+(X)$ such that $F(f)=\mu(f)-F^*(\mu)$. Then $F$ satisfies the downward continuity condition or equivalent to it, by Theorem \ref{compact}, for each $m\ge 0$, the set $\Sc_m=\{\mu\mid F^*(\mu)\le m\}$ is weak$^*$ compact.
\end{theorem}

{\bf Proof } Since $F$ satisfies the upward continuity condition we already know that
$$
F(f)=\sup \{\mu(f)-F^*(\mu)\mid \mu\in M_+(X)\}.
$$
The hypothesis says that this supremum is a maximum. As before we will use a compact metric space $K$ in which $X$ is dense.  We have already seen that when $f\in C(K)$:
$$
\sup \{\mu(f)-F^*(\mu)\mid \mu\in M_+(X)\}=\sup \{\mu(f)-F^*(\mu)\mid \mu\in M_+(K)\}.
$$
We will now suppose that $F$ does not satisfy the downward continuity condition and we hope that this leads to a contradiction. So we suppose that there is a sequence of functions $h_k\in C(K)$ and a $\delta>0$ such that $h_k\downarrow 0$ on $X$ but $F(h_k)\ge\delta$ for each $k$.  To keep the notation simple we use the same notation $h_k$ for the extension of $h_k$ to the compact space $K$.    Also on $K$ we have that $h_k$ is a nonincreasing sequence but  $g=\lim_kh_k$ is not necessarily $0$ on $K\setminus X$! Take now $m\ge \psi(\Vert h_1\Vert)$. To apply the sup-limsup theorem we use the following ``translation'':
\begin{eqnarray*}
E&=& \{\mu\in M_+(K)\mid F^*(\mu)\le m\}\\
G&=&\{\mu\in M_+(X)\mid F^*(\mu)\le m\}=\Sc_m\\
f_k&\colon& E\rightarrow \Rr\text{ is defined as } f_k(\mu)=\mu(h_k)-F^*(\mu).\\
\tilde{g}&\colon& E\rightarrow \Rr\text{ is defined as }\tilde{g}(\mu)=\lim_kf_k(\mu).
\end{eqnarray*}
The sequence $f_k$ is uniformly bounded since -- by the same argument as in Proposition \ref{bound} -- $E$ is a norm bounded set and for all $k$, $\Vert h_k\Vert\le\Vert h_1\Vert$. Also $\conv_\sigma\{f_k\mid k\ge 1\}\subset C(K)$.  It is easily seen that the hypothesis says that $G$ satisfies the peak condition.
By the monotone convergence theorem we have $\tilde{g}(\mu)=\mu(g)-F^*(\mu)$ for each $\mu\in E$. For each $k$ the set $L_k\subset E$ with $L_k=\{\mu\in E\mid f_k(\mu)\ge\delta\}$ is convex and compact since $f_k$ is concave and upper semi continuous for the weak$^*$ topology.  Furthermore each $L_k$ is nonempty since $F(h_k)\ge\delta$. The set $L=\cap_kL_k$ is nonempty since the sequence $(L_k)_k$ is also nonincreasing.  This means that there is $\nu\in E$ with $f_k(\nu)\ge\delta$ for each $k$. The monotone convergence theorem again says that $\nu(g)-F^*(\nu)\ge\delta>0$. This shows that
$$
\sup\{\limsup_kf_k(\alpha)\mid \alpha\in E \}\ge \delta>0.
$$
When we restrict to $X$ the situation is easier. On $G$ we have that $\mu(g)-F^*(\mu)\le 0$ and hence
$$
\sup\{\limsup_kf_k(\alpha)\mid \alpha\in G \}\le 0.
$$
This is a contradiction to the sup-limsup theorem. Until now we only showed that downward continuity  holds for sequences in $C(K)$.   We will now use the equivalence with weak$^*$ compactness as in Theorem \ref{compact} and show that for each $m$ the set $\Sc_m$ is uniformly tight or what is the same, weak$^*$ compact. We will actually show that
$$
\{\mu\in M_+(K)\mid F^*(\mu)\le m\}=\{\mu\in M_+(X)\mid F^*(\mu)\le m\}.
$$
Suppose on the contrary that $$\mu_0\in \{\mu\in M_+(K)\mid F^*(\mu)\le m\}\setminus\{\mu\in M_+(X)\mid F^*(\mu)\le m\}.$$ Then there is a compact set $L\subset K\setminus X$ with $\mu_0(L)>0$. Let now $h\in C(K)$ be such that $h=1$ on $L$ and $0\le h <1$ on $K\setminus L$.   Take now $k\ge 1$ so that $k\mu_0(L)-m\ge 1$. The sequence $k\,h^n$  converges downward to $0$ on $X$. According to the first part of the proof we have $F(k\,h^n)\rightarrow 0$. But this is in contradiction with $F(k\,h^n)\ge \mu_0(k\, h^n)-F^*(\mu_0)\ge k\mu_0(L)-m\ge 1$.\qed
\begin{remark}{\rm If $F$ is nondecreasing and satisfies the downward continuity property then for measures $\mu\in M(K)\setminus M(X)$ we have $F^*(\mu)=+\infty$.
}\end{remark}
\begin{remark}{\rm It is easily seen that the downward continuity property cannot be extended to a continuity for uniformly bounded pointwise converging sequences, the so called Lebesgue property. We just take a compact metric space $K$ and define $F(f)=\max_{x\in K}f(x)$ which can be written as $F(f)=\sup\{\mu(f)\mid \mu\in\Pc(K)\}$, i.e. the sup over all probabilities on $K$. Without giving details we mention that a necessary and sufficient condition for this Lebesgue property is that for each $m$, the set $\Sc_m$ is weakly compact, meaning it is compact for the topology $\sigma(M(X),M^*(X))$.
}\end{remark}

\section{The case of all continuous functions}

The result in this section is a generalisation of the result for coherent utility functions, \cite{C(X)}. It follows the same lines but the use of the monetary property will be replaced by another argument.  Also the theorem is formulated for convex functions and not just for positively homogeneous ones. For completeness we give the adapted proof.  The symbol $M^c(X)$ stands for measures on $X$ having a compact support, this is also the dual space of the space $C(X)$ with the usual topology of uniform convergence on compact sets.  We also observe that for $\mu\in M^c(X)$:
$$
\sup\{\mu(f)-F(f)\mid f\in C_b(X)\} = \sup\{\mu(f)-F(f)\mid f\in C(X)\}=F^*(\mu).
$$

\begin{theorem} Let $F\colon C(X)\rightarrow\Rr$ be a nondecreasing convex function with $F(0)=0$.  Then
\begin{enumerate}
\item $F$ satisfies an upward continuity property, i.e. for a sequence $(f_n)_n$ in $C(X)$ such that $f_n\uparrow f\in C(X)$  we have $F(f_n)\uparrow F(f)$.
\item There is  a compact set $L\subset X$ such  that the function $F$ factorises over the restriction map $C(X)\rightarrow C(L)$. In other words the conjugate function
\begin{eqnarray*}
F^*&\colon&M^c(X)\rightarrow \Rr_+\cup\{+\infty\}\\
& &\mu\rightarrow F^*(\mu)=\sup\{\mu(f)-F(f)\mid f\in C(X)\}
\end{eqnarray*}
can only take finite values on $M_+(L)$ and
$$
F(f)=\sup\{\mu(f)-F^*(\mu)\mid \mu\in M_+(L)\}.
$$
\item The function $F$ satisfies the downward continuity property since it factorises over the restriction map $C(X)\rightarrow C(L)$.
\end{enumerate}
\end{theorem}
{\bf Proof } Let a sequence $(f_n)_n$ in $C(X)$ be such that $f_n\uparrow f\in C(X)$ and $F(f)\ge 2\epsilon + F(f_n)$ where $\epsilon>0$. Since the function $\Rr\rightarrow \Rr; a\rightarrow F(f+a)$ is convex, it is continuous.  Therefore there is $\delta>0$ such that for $|a|\le\delta$, $|F(f+a)-F(f)|\le\epsilon$.  We replace the given sequence by $g=(f-\delta), g_n=\min(f-\delta,f_n)$.  We still have $F(g)\ge F(g_n)+\epsilon$.  But this time $g_n\uparrow g$ in a stationary way, i.e. for each $x\in X$ there is $n_0$ such that for $n\ge n_0$, $g_n(x)=g(x)$. The sets $G_n=\{x\mid g_n(x)<g(x)\}$ form a decreasing sequence of open sets and their closures satisfy $\overline{G_n}\subset \{x\mid f_n(x)\le f(x)-\delta\}$.  Their intersection is empty: $\cap_n \overline{G_n}=\emptyset$.  We now define $h_n=g+n(g-g_n)=(n+1)g-n\,g_n$ and observe that convexity implies
$$
F(g_n)+\epsilon\le F(g)\le \frac{n}{n+1} F(g_n)+\frac{1}{n+1}F(h_n).
$$
This simplifies to $F(h_n)\ge (n+1)\epsilon -F(g_n)\ge (n+1)\epsilon -F(g_1)$.  We now remark that $h_n\rightarrow g$ and that $h=\max (h_n; n\ge 1)$ is a real valued function.  We now show that $h$ is continuous.  To do this, take a sequence $x_m\rightarrow x$ forming the compact set $C=\{x_m; m\ge 1\}\cup\{x\}$.  The intersection $C\cap \cap_n \overline{G_n}=\emptyset$ and compactness implies that there is $n_0$ such that $\overline{G_{n_0}}\cap C=\emptyset$. This means that for all $x_m$ and for $x$ we have
\begin{eqnarray*}
h(x_m)&=&\max( h_n(x_m); n<n_0)\\
h(x)&=&\max( h_n(x); n<n_0).
\end{eqnarray*}
Since the maximum is taken over a finite set we obtain convergence $h(x_m)\rightarrow h(x)$.  Obviously $F(h)\ge F(h_n)\ge (n+1)\epsilon -F(g_1)$ which tends to infinity with $n$, violating that $F(h)$ must be finite.  This shows the first part of the theorem. For the second part  we start by defining the conjugate only using bounded functions. So for the moment we put
$$
F^*(\mu)=\sup\{\mu(f)-F(f)\mid f\in C_b(X)\}.
$$
We already know that $F^*$ can only take finite values for nonnegative measures.  Because of upward continuity we may recall from the previous section, that for $f\in C_b(X)$:
$$
F(f)=\sup \{\mu(f)-F^*(\mu)\mid \mu\in M_+(X), F^*(\mu)<\infty\}.
$$
The upward continuity and the monotone convergence theorem from measure theory then show that for elements $f\in C(X)$ that are bounded below, we get the same formula: $F(f)=\sup \{\mu(f)-F^*(\mu)\mid \mu\in M_+(X), F^*(\mu)<\infty\}.$  The next step is to show that everything happens on a compact subset of $X$. For $\mu\in M_+(X)$, let $S_\mu$ be the support of the measure $\mu$.  We claim that $\cup_{F^*(\mu)<\infty} S_\mu$ is a relatively compact set.  If not there would exist a sequence of measures $\mu_n$ with $F^*(\mu_n)<\infty$ as well as a sequence of different points $x_n\in S_{\mu_n}$ together with open balls $B(x_n,\eta_n)$ such that $\eta_n\rightarrow 0$ and such that for each subset $J\subset \Nr$, $\cup_{j\in J}\overline{B(x_j,\eta_j)}$ is closed and the closures $\overline{B(x_n,\eta_n)}$ are disjoint. For each $n$ we define $\psi_n(x)=\dist(x, X\setminus B(x_n,\eta_n))$. Observe that $\psi_n(x_n)>0$. Because $x_n\in S_{\mu_n}$ we have that $\mu_n(\psi_n)>0$.  Consequently there are numbers $a_n$ such that $a_n\mu_n(\psi_n)-F^*(\mu_n)\ge n$.  The function $h=\sum_n a_n\psi_n=\max(a_n\psi_n; n\ge 1)$ is still continuous (because of the choice of the numbers $\eta_n$). Now we have for each $n$:
$$
F(h)\ge F(a_n\psi_n)\ge a_n\mu_n(\psi_n)-F^*(\mu_n)\ge n
$$
This is a contradiction to the finiteness of the function $F$, showing that all the measures $\mu$ with $F^*(\mu)<\infty$ are supported by the same compact set, say $L$. This changes the dual representation of $F$ -- at least for functions bounded below -- into
$$
F(f)=\sup\{\mu(f\,\one_L)-F^*(\mu)\mid F^*(\mu)<\infty\}.
$$
We must still show that the same expression holds for functions not necessarily bounded below. Take an arbitrary function $f\in C(X)$ and for $m\le\inf_{x\in L}f(x)=m_0$, observe that all the expressions $F(f\vee m)$ are the same, say $\alpha$.  Suppose now that $F(f) <\alpha -\delta$ where $\delta>0$. The function $f'=2( f\vee m_0) - f = f\vee m_0 +(f\vee m_0-f)$ is bounded below and on $L$ it equals $f$. Convexity then gives
$$
\alpha=F(f\vee m_0)\le {1\over 2} F(f') + {1\over 2} F(f) \le {1\over 2} (\alpha +\alpha -\delta)=\alpha-\frac{\delta}{2},
$$
a contradiction.  Hence for all $f$ we get: 
$$
F(f)=\sup\{\mu(f\,\one_L)-F^*(\mu)\mid \mu\in M_+(L)\},
$$
and $F$ factorises over the restriction $C(X)\rightarrow C(L)$.  The last statement in the theorem is now straightforward. \qed

\end{document}